\documentclass[12pt]{article}
\usepackage{amsmath}
\usepackage{graphicx,psfrag,epsf}
\usepackage{enumerate}
\usepackage{natbib}
\usepackage{amsmath,amsthm,amsfonts,epsfig,amssymb,natbib,eucal,eufrak}
\usepackage{booktabs}
\usepackage{multirow}
\usepackage{siunitx}
\usepackage{qtree}
\usepackage{color}
\usepackage{setspace}

\usepackage{float}
\restylefloat{table}

\newtheorem{lem}{Lemma}

\newcommand{\blind}{0}

\begin{document}

\def\spacingset#1{\renewcommand{\baselinestretch}%
{#1}\small\normalsize} \spacingset{1}


\if0\blind
{
  \title{\bf Proportional Closeness Estimation of Probability of Contamination Under Group Testing}
  \author{Yaakov Malinovsky\thanks{Corresponding author}\\
    Department of Mathematics and Statistics\\ University of Maryland, Baltimore County, Baltimore, MD 21250, USA\\
    and \\
    Shelemyahu Zacks\\
    Department of Mathematical Sciences\\ Binghamton University,
Binghamton, New York, USA
    }

  \maketitle
} \fi

\if1\blind
{
  \bigskip
  \bigskip
  \bigskip
  \begin{center}
    {\LARGE\bf Title}
\end{center}
  \medskip
} \fi

\doublespacing

\bigskip
\begin{abstract}
The paper is focused on the problem of estimating the probability $p$ of individual contaminated sample, under group testing. The precision of the estimator is given by the probability of proportional closeness, a concept defined in the Introduction. Two-stage and sequential sampling procedures are characterized. An adaptive procedure is examined.
\end{abstract}

\noindent%
\vfill


\section{Introduction}
\label{se:I}
Let $p$ be the probability of contamination of an individual sample. We assume that all individual samples have the same $p$ and are independent. Thus, the samples represent Bernoulli trials for estimating $p$. In group testing the samples of several individuals, $k$ say, are mixed and tested together for contamination. The result of the test is binary, $1$ if the group material is contaminated and $0$ otherwise.
Thus, the number of contaminated groups, among $n$ independent ones, has a binomial distribution with parameters $\theta_k=1-(1-p)^k$ and $n$.
The goal is to estimate the individual probability $p$ from the data of group testing.

The earliest study using group testing, according to \cite{Boswell1996}, is ascribed to \cite{Watson1936}. Early studies of group testing deal with the estimation of the rate of disease transmission, from insects to plants \citep{GG1960, CR1962, Thompson1962}. These studies deal with the maximum likelihood estimation (MLE) of $p$ based on a fixed number of groups. Bias correction is studied, as well as, the optimal sample size, as a function of $p$, to minimize the mean squared error (MSE) of the estimator. Confidence intervals for $p$ are given too. Progress in the estimation based on fixed sample size was achieved by \cite{Bur1987}, who proposed an estimator with bias proportional to the reciprocal of the squared sample size. This estimator has an MSE smaller than that of the MLE. Comparison of different estimators based on fixed sample size can be found in \cite{HW2009}. Under fixed sample size there is no unbiased estimator of $p$ \citep{B1954}. \cite{Hall1963} studied estimation under inverse binomial sampling. \cite{Haber2017} constructed an unbiased estimator, based on inverse binomial sampling, with a stopping rule based on reaching a given number of negative (uncontaminated) groups.   Sequential sampling designs were also proposed by \cite{Kerr1971,KU2006,PT2011,H2013}. Sequential methods depend on the type of stopping rules applied.

In the present paper we study the properties of two-stage and sequential sampling with stopping rules based on prescribed probability of proportional closeness estimators.
Such stopping rules were not considered before in the context of group testing.
We also assume that the testing for contamination is error free. Prescribed probability of proportional closeness (PPPC) estimation of a parameter $\theta$, by an estimator $\hat{\theta}$ should satisfy the following probability requirement.
For a given $0<\alpha,\,\gamma<1$,

 \begin{equation*}
 \displaystyle
 P_{\theta}\left(-\gamma \theta< {\hat{\theta}-\theta}<\gamma \theta \right)\geq 1-\alpha,\,\,\, \text{for all}\,\,\,\,0<\theta<1.
 \end{equation*}
 In estimating the probability $p$, we wish to design a two-stage or a sequential procedure with a stopping variable $N$, so that the estimator at stopping will satisfy,
for a given $0<\alpha,\,\gamma<1$,
 \begin{equation}
 \label{eq:c2}
 \displaystyle
 P_{p}\left(-\gamma<\frac{\hat{p}_{N}-p} {p}<\gamma \right)\geq 1-\alpha,\,\,\, \text{for all}\,\,\,\, 0<p<1.
 \end{equation}

The PPPC criterion for proportional closeness was introduced by \cite{EL1964} (p. 339), and was applied by
\cite{Z1966} for sequential estimation of the mean of a log-normal distribution. This criterion was used also later by
\cite{N1969}, \cite{WF1983}, and others. Proportional closeness estimation of quantities is being used widely in physical sciences and engineering,
where accuracy is measured as proportion of the true value of the measured quantity. Notice that implied interval $\displaystyle \left((1-\gamma)p,\,(1+\gamma)p\right)$ is a prediction interval for the estimator $\displaystyle \hat{p}_{N}$, not a confidence interval
for $p$. In order to obtain an estimator inside this prediction interval, since the estimator is consistent, one needs to determine the sample size
$N$ to satisfy the probability statement in \eqref{eq:c2}. The sample size for this criterion is function of $p$. Since the value of
$p$ is unknown, we need two-stage or sequential sampling, in order to obtain closeness probability approximately equal to the prescribed $\displaystyle 1-\alpha$.

An alternative approach could be to apply the dual proportional closeness confidence interval
$\displaystyle {
\left(\hat{p}_{N}/(1+\gamma),\,\min\left\{1,\,\hat{p}_{N}/(1-\gamma)\right\}\right)
}$, and determine the sample size $N$ so that the coverage probability of $p$ will be approximately  $\displaystyle 1-\alpha$.
This approach is similar to the one used by \cite{MB2014, MB2015, MB2016} and \cite{DM2015}, called
 \textquotedblleft fixed accuracy confidence interval \textquotedblright\ which is
 $\displaystyle {\left(\hat{p}_{N}/\delta,\, \min\left\{1,\,\hat{p}_{N}\delta\right\} \right)}$
 for some $\delta>1$. Notice that if $\displaystyle \delta=1+\gamma$ and $\displaystyle \gamma$ is small, then the two types
 of confidence intervals are almost the same. For the fixed accuracy confidence intervals see also the paper of \cite{MZ2016}.

In the present paper we are studying the problem of estimating the prevalence $p$, when $p$ is small.
It is motivated by a need to estimate prevalence of infectious diseases in the medical studies.
For example, the group testing is used to estimate prevalence of HIV \citep{P2005} (overall reported rate in the study was 0.1\%) and Hepatitis B virus \citep{S2013} (overall reported rate in the study was 0.01\%).

All our examples are for $p\leq 0.5$.
When $p=0.5$, the required sample, for $\alpha=0.05,\, \gamma=0.1$, is $385$ (see Table\ref{t:1}).
In this case the probability that $\displaystyle \hat{p}_{N}>0.9$ is approximately $\displaystyle \Phi\left(1.8/\sqrt{385}-\sqrt{385}\right)\approx 0$. Accordingly, in the following we consider the proportional accuracy interval  $\displaystyle {
\left(\hat{p}_{N}/(1+\gamma),\,\hat{p}_{N}/(1-\gamma)\right)
}$
with $\gamma=0.1$.

We show first how large should a sample of individual observations be ($k=1$), in order to satisfy the prescribed proportional coverage probability (PPCP) requirement. We then study the group testing properties.

\section{Large Sample Approximation}
\subsection{Individual Testing}
Let $\displaystyle\left\{J_i,\,i=1,\ldots,n\right\}$
be i.i.d. random binary variables, with $\displaystyle P\left(J_i=1\right)=p=1-P\left(J_i=0\right).$
For a random sample of $n$ Bernoulli trials, the minimal sufficient statistic is $$\displaystyle \widehat{p}_{n}=\frac{1}{n} \sum_{i=1}^{n}J_i=_{d}\,\,\frac{1} {n}Bin\left(n,p\right).$$
By Central Limit Theorem (CLT) $\displaystyle\sqrt{n}\left(\widehat{p}_{n}-p\right)\rightarrow_{d}N\left(0,\, pq \right),$ where $q=1-p$.
As will be shown, to satisfy the proportional coverage condition \eqref{eq:c2}, when $p$ is small, large samples are required. Accordingly, we apply the large sample normal approximation. For large samples, we have,
\begin{equation}
\label{eq:ls}
P_{p}\left(-\gamma<\frac{\hat{p}_{n}-p}{p}<\gamma\right)\cong 2 \Phi\left(\gamma\sqrt{\frac{np} {q}}\right)-1,
\end{equation}
where $\displaystyle \Phi(\cdot)$ denotes the standard normal distibution.
Thus, from the large samples approximation \eqref{eq:ls}, it follows that the coverage probability \eqref{eq:c2} will be approximately satisfied if $\displaystyle n>n^{*}\left(p,\,\gamma\right)$, where
\begin{equation*}
\displaystyle
{ n^{*}\left(p,\,\gamma\right)=\chi^{2}_{1-\alpha}\frac{q} {p\gamma^2}.}
\end{equation*}
$\displaystyle \chi^{2}_{1-\alpha}$ denotes the $\displaystyle (1-\alpha)$-quantile of the chi-squared distribution with 1 degree of freedom.
In Table \ref{t:1} we display the values of $\displaystyle n^{*}\left(p,\,\gamma\right)$, when $\displaystyle \alpha=0.05$ and $\gamma=0.1$.

\begin{table}[H]
\caption{Required Sample Size (Large Sample Approximation)}
\label{t:1}
\begin{center}
\begin{tabular} {l|llllllllllll}
 \hline
 $p$ &0.01&0.05& 0.1& 0.2& 0.3& 0.4& 0.5& 0.6& 0.7& 0.8& 0.9\\
 $n^{*}\left(p,\,\gamma\right)$ &38031&7299&3458&1537&897&577&385&257&165&97&43\\
 \hline
 \end{tabular}
 \end{center}
 \end{table}

 \subsection{Group Testing}
 In group testing (GT) the material of $k$ individuals is mixed together into one batch, and one test is performed on this batch for contamination. This GT can be repeated independently several times in order to estimate the probability of batch contamination $\displaystyle \theta_k$, and from this to obtain and estimator of $p$. We can model this as Bernoulli trials with batches of size $k$. Let $X=1$ be a random variable signifying that the batch is contaminated, and $X=0$ if the batch is not contaminated. Let
 $\theta_{k}=P\left(X=1\right)$. Notice that $\theta_{k}=1-(1-p)^{k}$.
Let $S_n=\sum_{i=1}^{n}X_i$.
Then, the distribution of $S_{n}$ is $\displaystyle Bin(n,\theta_k)$. The
maximum likelihood estimator (MLE) of $\theta_k$ is
$$\displaystyle \widehat{\theta}_{n,\,k}=\frac{S_n}{n}=\overline{X}_{n,\,k}.$$
The MLE of $p$ after testing $n$ groups is
\begin{equation*}
\displaystyle \widehat{p}_{n,\,k}=1-\left(1-\overline{X}_{n,\,k}\right)^{1/k}.
\end{equation*}
Notice that a.s. $\displaystyle \lim_{n\rightarrow\infty} \widehat{p}_{n,\,k}=p$.
Furthermore,
\begin{equation*}
E \left(\overline{X}_{n,\,k}\right)=1-q^{k},\,\,\,\,
Var\left(\overline{X}_{n,\,k}\right)=\frac{q^{k}-q^{2k}}{n}.
\end{equation*}

Define the function $$\displaystyle g(x)=1-(1-x)^{1/k},\,\,0<x<1.$$
Notice that $\displaystyle p=g(\theta_{k})$, and $\displaystyle
\widehat{p}_{n,\,k}=g\left(\overline{X}_{n,\,k}\right)$.
According to the delta method, for large values of $n$, we get
\begin{equation}
\label{eq:e}
\displaystyle
E\left(\widehat{p}_{n,\,k}\right)
=p+\frac{e_k\left(q\right)}{n}
+o\left(\frac{1}{n}\right),
\end{equation}
where $\displaystyle e_k(q)=\frac{k-1}{{2k^2}}q\left(q^{-k}-1\right)$,
and
\begin{equation}
\label{eq:v}
\displaystyle Var\left(\widehat{p}_{n,\,k}\right)=
\frac{v_k\left(q\right)} {n}+0\left(\frac{1} {n^2}\right),
\end{equation}
where $\displaystyle  v_k\left(q\right)=\frac{q^2(q^{-k}-1)} {k^2}$.

Finally, the asymptotic distribution of $\displaystyle \widehat{p}_{n,\,k}$ is normal. Thus
\begin{equation}
\label{eq:lsa}
\displaystyle
P_{p}\left(-\gamma<\frac{\widehat{p}_{n,\,k}-p} {p}<\gamma\right)\cong
2\Phi\left(\sqrt{n}\frac{p\gamma}{\sqrt{v_k\left(q\right)}}\right)-1.
\end{equation}

It follows from \eqref{eq:lsa} that a large sample approximation for the $1-\alpha$ coverage probability is attained if $\displaystyle n>n^{*}_{G}\left(p ,k\right),$ where
\begin{equation}
\label{eq:nG}
\displaystyle n^{*}_{G}\left(p ,k\right)=\frac{\chi^{2}_{1-\alpha}}{\gamma^2}\frac{v_k\left(q\right)}{p^2 }=\frac{\chi^{2}_{1-\alpha}}{(\gamma k)^2}\psi\left(\theta_k\right),\,\,
\end{equation}
and
\begin{equation*}
\psi\left(\theta_k\right)=\frac{\theta_k \left(1-\theta_k\right)^{2/k}}{\left(1-\theta_k\right) \left(1-\left(1-\theta_k\right)^{1/k}\right)^2}.
\end{equation*}
In Table 2 we present a few values of $\displaystyle n^{*}_{G}\left(p ,k\right)$ for various values of $p$ and $k$, when $\alpha=0.05, \gamma=0.1$.

\begin{table} [H]
\caption{ Minimal Number of Groups Required}
\label{t:2}
\begin{center}
\begin{tabular} {l|llllllll}
  \hline
  $k$ &1 &50 &100 &158 & 159 & 160 & 170 & 200 \\
  $p=0.01$ &38030.44 & 983.23 & 652.10 & 587.242 &{\bf 587.24} & 587.27 & 588.98 &608.41\\
  \hline

  $k$ &1 &10 &20 &30 &31 &32 & 33 &40 \\
  $p=0.05$ &7298.80 &929.39 &620.41 &563.80 & {\bf 563.39} &563.68 &564.64 &587.76\\
  \hline

  $k$ &1 &10 &11 &12 &13 &14 &15 &16\\
  $p=0.10$ &3457.30 &581.20 &562.30 &549.00 &540.20 & 535.20 &{\bf 533.40} &534.40\\
  \hline
  $k$ &1 &6 &7 &8 &9 &10 &11 &12 \\
  $p=0.20$ &1536.60 &480.60 &{\bf 472.70} &476.40 &489.50 &511.00 &540.50 &578.40 \\
  \hline
  $k$ &1 &2 &3 &4 &5 &6 &7 &8\\
  $p=0.30$ &896.30 &544.20 &445.10 &{\bf 413.70} &414.10 &435.70 &475.60 &534.20 \\
  \hline
  $k$ &1 &2 &3 &4 &5 &6 &7 \\
  $p=0.40$ & 576.20&384.10 &{\bf 348.60} &362.80 &410.00 &490.60 &612.50 \\
  \hline
  $k$ &1 &2 &3 &4 &5 &6 &7 \\
  $p=0.50$ & 384.10&{\bf 288.1} &298.80 &360.10 &476.30 &672.30 &995.60 \\
  \hline
  \end{tabular}
  \end{center}
  \end{table}
  \smallskip
We see in Table \ref{t:2} that, if we know the value of $p$ we can determine the optimal $k$, which minimizes $n^{*}_{G}\left(p ,k\right)$. Moreover, if
the test of a group costs the same as that of an individual sample, the GT achieves the same PPCP with significantly smaller cost, provided the cost of collecting the individual samples is relative small.

\section{Sequential Procedures }
In a sequential procedure, one takes first m groups of size k, estimates
$\widehat {p}_{m,\,k}$ and makes a decision whether to stop testing, or continue by taking one group at a time. After testing $n\geq{m}$ groups, one needs a stopping rule, based on the statistic $\overline{X}_{n,\,k}$. We start first by defining a stopping rule for a fixed $k$, and studying it's properties. We suggest the following stopping variable, based on \eqref{eq:nG}, namely

\begin{equation}
\label{eq:i}
\displaystyle
N=\min\left\{n\geq m:\,n>\frac{\chi^{2}_{1-\alpha}}{(\gamma k)^2}\psi\left(\overline{X}_{n,\,k}\right)\right\}.
\end{equation}

Recall that $\overline{X}_{n,\,k}$ is an unbiased estimator of
$\theta_{k}$. We will obtain the asymptotic distribution of the stopping variable $N$.

By the delta method we derive the asymptotic approximation of
$ E\left(\psi\left(\overline{X}_{n,\,k}\right)\right)$, which is
\begin{equation*}
\displaystyle
E\left(\psi\left(\overline{X}_{n,\,k}\right)\right)\approx
\psi\left(\theta_k\right)+\frac{1}{2 n} \psi^{''}\left(\theta_k\right)
\theta_k\left(1-\theta_k\right).
\end{equation*}
Similarly, the asymptotic variance of $\psi\left(\overline{X}_{n,\,k}\right)$ is
\begin{equation*}
\displaystyle
AV\left(\psi\left(\overline{X}_{n,\,k}\right)\right)\approx
\frac{1}{n}\left(\psi^{'}\left(\theta_k\right)\right)^2
\theta_k\left(1-\theta_k\right).
\end{equation*}

Since the distribution of $\psi\left(\overline{X}_{n,\,k}\right)$ is normal, the distribution of $N$ is obtained in the following way. Let
\begin{equation*}
\displaystyle
\zeta_k=\frac{\chi^2_{1-\alpha} }{\left(\gamma k\right)^2}.
\end{equation*}
We obtain that
\begin{equation}
\label{eq:17}
P\left(N>n\right)=P\left(\frac{n}{\zeta_k}<\psi\left(\overline{X}_{n,\,k}
\right)\right)\approx 1-\Phi\left(\frac{n-\zeta_k \psi\left(\theta_k\right)}{\zeta_k \sqrt {AV\left(\psi\left(\overline{X}_{n,\,k}\right)\right)}}\right).
\end{equation}
Furthermore,
\begin{align}
\label{eq:m}
&
P\left(N=m_k\right)=1-P\left(N>m_k\right),\,\,\,\text{for}\,\,\, n\geq m_k,\\ \nonumber
&
P\left(N=n\right)=P\left(N>n-1\right)-P\left(N>n\right).
\end{align}

The moments of $N$ can then be computed from the probability mass function (p.m.f) \eqref{eq:m}.
The maximum likelihood estimator (MLE) of $p$ at stopping is
\begin{equation*}
\displaystyle
\widehat{p}_{N,\,k}=1-\left(1-\overline{X}_{N,\,k}\right)^{1/k}.
\end{equation*}
Using \eqref{eq:e} and \eqref{eq:v} we obtain the asymptotic approximations:
\begin{equation*}
\displaystyle
E\left(\widehat{p}_{N,\,k}\right)=p+e_{k} E\left(\frac{1}{N}\right),
\end{equation*}
and
\begin{equation*}
\displaystyle
Var\left(\widehat{p}_{N,\,k}\right)=e_k^2 Var\left(\frac{1}{N}\right)+v_kE\left(\frac{1}{N}\right).
\end{equation*}

The coverage probability of the proportional closeness interval is
\begin{equation*}
\displaystyle
CP=P\left(-\gamma<\frac{\widehat{p}_{N,\,k}-p}{p}< \gamma\right)
\approx 2E\left(\Phi\left(\gamma p \sqrt{\frac{N}{v_k}}\right)\right)-1.
\end{equation*}

Large sample approximations and simulated estimates of the characteristics of the sequential procedure are displayed in the following table.

\begin{table}[H]
\label{t:3}
\caption{Estimates (A= analytic approximation, S= simulation based on 1000 replicates) under optimal conditions of the sequential procedure, with $\alpha=0.05, \gamma=0.1$ }
\begin{center}
\begin{tabular}{l|lllllllllllllll}
  \hline
   & $p$ & $k$ & $m_k$ & $E\left(N\right)$ & $\sigma\left(N\right)$ & $E\left(\widehat{p}_{N,\,k}\right)$ & $\sigma\left(\widehat{p}_{N,\,k}\right)$ & CP \\
  \hline
  A&0.5 & 2   & 250 & 288.23 & 19.06 & 0.5007& 0.0256& 0.9494 \\
  S&0.5 & 2   & 250 &299.83&   18.68&   0.5014&  0.0260&    0.9430 \\
  A&0.4 & 3   & 320 & 348.93 & 13.64 &0.4007& 0.0204& 0.9499\\
  S&0.4 & 3   & 320  &348.86&   13.41&    0.4018&    0.0199&    0.9500  \\
  A&0.3 & 4   & 390 & 414.18 & 12.32 &0.3005& 0.0153& 0.9500\\
  S&0.3 & 4   & 390  &    414.50&   12.33&    0.3011&    0.0153&    0.9470    \\
  A&0.2 & 7   & 450 & 473.14 & 6.89  &0.2004& 0.0102& 0.9501\\
  S&0.2 & 7   & 450  &    473.50&    6.60&    0.2009&    0.0099&    0.9600    \\
  A&0.1 & 15  &510  & 533.87 & 3.37 & 0.1002& 0.0051& 0.9501\\
  S&0.1 & 15  &510   &    534.77&   3.58&     0.1001&    0.0053&    0.9390    \\
  A&0.05 &31  &550 &  563.89& 1.63  & 0.0501& 0.0025& 0.9501\\
  S&0.05 &31  &550   &    564.63&    2.00&    0.0502&    0.0026&    0.9530    \\
  A&0.01 &159 &585 & 587.88 & 1.33& 0.0100&0.0005& 0.9501 \\
  S&0.01 &159 &585   &    588.75&    1.27&    0.0100&    0.0005&   0.9470    \\
   \hline
\end{tabular}
\end{center}
\end{table}

\section{Two-Stage Procedure }
In a two-stage sampling procedure, we chose first $m$ for the number of groups to sample in stage 1, and the size $k_1$ of these groups. We then compute the estimator $\overline{X}_{m,\,k_1}$ and the required total number of groups $N_{m,\,k_1}$, which is like in \eqref{eq:i}
\begin{equation}
\label{eq:25}
N_{m,\,k_1}=\zeta_{k_1}\psi\left(\overline{X}_{m,\,k_1}\right).
\end{equation}

If $\displaystyle m>N_{m,\,k_1}$ sampling is stopped after stage 1; otherwise,
$$\displaystyle M=\lfloor N_{m,\,k_1}\rfloor+1-m$$ new groups are independently tested in stage 2 of sampling,
where $\displaystyle \lfloor x\rfloor$ for $x>0$ is defined as the largest integer which is smaller or equal to $x$.

The size of each group in stage 2 is $k_2\left(\widehat{p}_1\right)$ (determined according to Table 2), where $\widehat{p}_1$ is the estimator of $p$ after stage 1. The required total number of groups is
\begin{equation}
\displaystyle
N_{2}=m 1_{\left\{m>N_{m,\,k_1}\right\}}+\left(M+m\right) 1_{\left\{m\leq N_{m,\,k_1}\right\}}.
\end{equation}

\subsection{MLE of $p$}

To simplify notation, Let $\overline{X}_1$ and $\overline{X}_2$ denote the estimators from stage 1 and 2 respectively. Similarly, let $k_1$ and $k_2$ denote the k-values used in stage 1 and stage 2.
If sample stops after stage 1 then the MLE of $p$ is
\begin{equation*}
\displaystyle
\widehat{p}_{MLE}=1-\left(1-{\overline{X}_1}\right)^{1/k_1}.
\end{equation*}
On the other hand if sample includs also stage 2, then the MLE of $p$ can be calculated as follows.

If $\displaystyle k_1=k_2=k$ then the MLE of $p$ is
\begin{equation*}
\displaystyle
\widehat{p}_{MLE}=1-\left(1-\overline{\overline{X}}\right)^{1/k},
\end{equation*}
where
\begin{equation*}
\overline{\overline{X}}=\frac{m \overline{X}_1+M \overline{X}_2}{N_2}.
\end{equation*}
On the other hand, if $k_1\neq k_2$ then the MLE is the maximizer of the likelihood function
\begin{equation*}
\displaystyle
L\left({p}\right)= \left(1-\left(1-p\right)^{k_1}\right)^{S_1}
\left(1-\left(1-p\right)^{k_2}\right)^{S_2} \left(1-p\right)^D,
\end{equation*}
where $S_1=m \overline{X}_1$, $S_2=M \overline{X}_2$ and
$D=k_1 \left(m-S_1\right)+k_2 \left(M-S_2\right)$.

Let $l\left(p\right)=logL\left(p\right)$.
In the following lemma the MLE is obtained by finding the zero of the score function $\displaystyle l^{'}\left(p\right)$, which is polynomial in $q=1-p$.
This lemma is very close to Theorem 1 by \cite{Hardwick1998}.

\begin{lem}
MLE of $q$ is unique root in $[0,1]$ of the equation
\begin{align*}
&
Aq^{k_1+k2}+Bq^{k_1}+Cq^{k2}+D=0,
\end{align*}
where $D=k_1(m-S_1)+k_2(M-S_2),\, C=-D-k_2S_2,\, B=-D-k_1S_1,\, A=D+k_1S_1+k_2S_2$.
\end{lem}

In the following table we evaluate the performance of two-stage procedure under limited information ($m=100,\, k_1=2$) compared with the optimal
$k_1$ and $m$ (from Table \ref{t:2}). Estimation is based on simulation with 1000 replicas, and $\alpha=0.05$ and $\gamma=0.1$.
\begin{table}[H]
\caption{Two-Stage Procedure, with $\alpha=0.05, \gamma=0.1$ }
\label{t:4}
\begin{center}
\begin{tabular}{lllllllllllllll}
  \hline
  $p$ & $k_1$ & $m$ &$E\left(N_{2}\right)$ & $\sigma\left(N_{2}\right)$& $E\left(\widehat{p}_{N_{2}}\right)$  & $\sigma\left(\widehat{p}_{N_{2}}\right)$ & CP\\
  \hline
  0.5&2&100& 289.094&   34.746&0.5030&    0.0258&    0.946&\\
  0.5&2&200&289.249&23.546& 0.5027&0.0251&0.952&\\
  \hline
  0.4&2&100&  389.975&46.811& 0.4021&0.0213&0.934\\
  0.4&3&200& 350.810&18.752&0.4009&0.0209&0.946\\
  \hline
  0.3&2&100&549.122& 77.186&  0.3012&0.0140& 0.963\\
  0.3&4&300&414.733&14.823&0.3010&0.0148&0.957\\
  \hline
  0.2&2&100& 884.471& 153.775&0.2011 &    0.0082& 0.978\\
  0.2&7&400&473.873&7.281&0.2002&0.0101&0.960\\
  \hline
  0.1&2&100&1929.4&487.1&0.1002& 0.0028 &0.998\\
  0.1&15&500&534.794&3.575&0.1003&0.0050&0.953\\
  \hline
  0.05&2&100& 4164.3& 1646.0& 0.0500&    0.0010&    1.000\\
  0.05&31&500&564.676&2.081&0.0502&0.0025&0.954\\
  \hline

\end{tabular}
\end{center}
\end{table}

In table \ref{t:4} we realize how important it is to start with the optimal $k_1$.
However, it cannot be done if $p$ is unknown. This problem will be addresed later.
\bigskip

\noindent
{\it The Fisher Information in the Two-Stage Procedure}\\
The asymptotic distribution of the MLE is normal with mean $p$ and variance which is equal to the inverse of the Fisher Information (FI) function. We derive now the FI function in a two-stage case. It helps to validate the results in the Table \ref{t:4}.
The expected value of the score function is 0. Indeed
\begin{equation}
\displaystyle
E\left( l^{'}\left(p\right)|\overline{X}_1\right)= \frac{S_1 k_1 \left(1-p\right)^{k_1 -1}} {1-\left(1-p\right)^{k_1}}-\frac{\left(m-S_1\right) k_1}{1-p}.
\end{equation}
From this we obtain that the expected score is 0.
The FI function is then defined
\begin{equation}
\label{eq:vv}
FI\left(p\right)=V\left(l^{'}\left(p\right)\right).
\end{equation}
Moreover,
\begin{equation}
\label{eq:ev}
E\left(V\left(l^{'}\left(p\right)|\overline{X}_1\right)\right)=
\frac{k_2^2 \left(1-p\right)^{k_2-2}E\left(M\right)}{1-\left(1-p\right)^{k_2}}.
\end{equation}
Also,
\begin{equation}
\label{eq:ve}
V\left(E\left(l^{'}\left(p\right)|\overline{X}_1\right)\right)=
\frac{m \left(k_1\right)^2 \left(1-p\right)^{k_1-2}}{1-\left(1-p\right)^{k_1}}
\end{equation}
The sum of Eq. \eqref{eq:ev}-\eqref{eq:ve} yields an explicit formula for Eq. \eqref{eq:vv}. In the following table we present the exact values of FI and the corresponding $\sigma\left(\widehat{p}\right)$
\begin{table}[H]
\caption{Exact Fisher Information and Asymptotic STD}
\begin{center}
\begin{tabular}{llllll}
\hline
$p$&$m$&$k_1$&$k_2$&FI&$\sigma\left(\widehat{p}\right)$\\
\hline
0.5&200&2&3&1522.37&0.02563\\
0.4&200&3&3&2410.58&0.02037\\
0.3&300&4&4&4273.46&0.01529\\
0.2&400&7&7&9615.70&0.01098\\
0.1&500&15&16&38446.0&0.00510\\
0.05&500&31&31&152783&0.00255\\
\hline
\end{tabular}
\end{center}
\end{table}

\subsection{Linear Estimator of Two MLE's}
The MLE estimators of $p$, based on the results of stage 1 and stage  2 are $\widehat{p}_1=1-\left(1-\overline{X}_1\right)^{1/k_1}$ and
$\widehat{p}_2=1-\left(1-\overline{X}_2\right)^{1/k_2}$ . These two estimators are not independent. We consider here a linear combination of these two MLE's, as an alternative for the grand MLE discussed in the previous section, namely
\begin{equation}
\displaystyle
\widehat{p}_{N_2}=\frac{\left(m \widehat{p}_1+M \widehat{p}_2\right)}
{N_2}.
\end{equation}
We wish to investigate how good is this alternative estimator.
\bigskip

{\it{The distribution of $N_2$}}
\smallskip

 We approximate the distribution of $N_2$ by the distribution of $N_{m,\,k_1}$. This approximation will be compared to results of simulations, in order to assess the goodness of the approximation.
As in Eq.\eqref{eq:17}, the distribution of $N_{m,\,k_1}$ depends on
$\overline{X}_1$ and is given by formula \eqref{eq:25}, i.e.,
\begin{equation}
\displaystyle
p\left(N_2\leq n\right)\approx \Phi\left(\frac{n-\zeta_{k_1} \psi\left(
\theta_{k_1}\right)}{\zeta_{k_1} \sqrt{AV\left(\psi\left(\overline{X}_1\right)\right)}}\right).
\end{equation}

We derive the expected value and variance of $\widehat{p}_{N_2}$ in Appendices \ref{A:a} and \ref{A:b}.

The coverage probability CP is
\begin{equation*}
\displaystyle
CP=P\left(p \left(1-\gamma\right)<\widehat{p}_{N_2}<p \left(1+\gamma\right)\right).
\end{equation*}
Let $B\left(\widehat{p}_{N_2}\right)=E\left(\widehat{p}_{N_2}\right)-p$.Then, since the asymptotic distribution of $\widehat{p}_{N_2}$ is normal,the coverage probability is approximately
\begin{equation*}
\displaystyle
CP=\Phi\left(\frac{\gamma p-B\left(\widehat{p}_{N_2}\right)}{SE\left(\widehat{p}_{N_2}\right)}\right)
+\Phi\left(\frac{\gamma p+B\left(\widehat{p}_{N_2}\right)}{SE\left(\widehat{p}_{N_2}\right)}\right)
-1,
\end{equation*}
where $SE\left(\widehat{p}_{N_2}\right)=\left(V\left(\widehat{p}_{N_2}
\right)\right)^{1/2}$.
In the following table we present the exact and simulated functionals of the two-stage sampling and the linear combination of MLE's.
\begin{table}[H]
\caption{Exact and Simulated values of Two-Stage Functionals With Linear Combination of MLE's}
\begin{center}
\begin{tabular}{l|llllllllllll}
\hline
$$&  $p$ & $k_1$ & $m$ & $E\left(N_{2}\right)$ & $\sigma\left( {N_{2}}\right)$ & $E\left(\widehat{p}_{N_2}\right)$ & $SE\left(\widehat{p}_{N_2}\right)$ & $CP$\\
\hline
E&0.5&2&200&288.61&23.526&0.498&0.02868&0.919\\
S&0.5&2&200&289.68&23.639&0.501&0.03400&0.837\\
\hline
E&0.4&3&200&349.08&18.370&0.400&0.02428&0.901\\
S&0.4&3&200&349.74&18.282&0.402&0.03008&0.826\\
\hline
E&0.3&4&300&414.21&14.905&0.300&0.01880&0.889\\
S&0.3&4&300&414.73&14.224&0.301&0.01913&0.893\\
\hline
E&0.2&7&400&473.19&7.493&0.200&0.01215&0.891\\
S&0.2&7&400&473.27&7.548&0.201&0.01230&0.891\\
\hline
E&0.1&15&500&533.88&3.481&0.100&0.00605&0.902\\
S&0.1&7&500&534.25&3.748&0.100&0.00559&0.928\\
\hline
E&0.05&31&500&563.89&1.726&0.050&0.00293&0.912\\
S&0.05&31&500&564.32&2.054&0.050&0.00284&0.919\\
\hline
\end{tabular}
\end{center}
\end{table}

\section{Adaptive Designs}
In practice, the true value of $p$ is unknown. We have seen that the optimal size of groups $k$ depends on $p$. We therefore suggest to start with a first pilot sample of $m_0=100$ groups of size $k_0=2$. The value of $p$ is then estimated by the corresponding MLE $\hat{p}_0$. The value of this estimator can be used as a required parameter for determinig $k=kp(\hat{p}_0)$ and the corresponding $m=mp(\hat{p}_0)$. In the following table we present simulation estimates of such an adaptive sequential procedure. In these simulations the number of independent runs is 1000. Also $\alpha=0.05$ and $\gamma=0.1$.

\begin{table}[H]
\caption{Simulated estimates of the adaptive sequential procedure }
\begin{center}
\begin{tabular}{lllllllllllllll}
  \hline
  $p$ & $k_0$ & $m_0$ & $E\left(N_{3}\right)$ & $\sigma\left(N_{3}\right)$& $E\left(\widehat{p}_{N_{3}}\right)$  & $\sigma\left(\widehat{p}_{N_{3}}\right)$ & CP\\
  \hline
  0.5&2&100&403.271&20.291&0.5018&0.022&0.977\\
  0.4&2&100&462.828&18.431&0.4019&0.018&0.976\\
  0.3&2&100&524.398&15.366&0.3014&0.014&0.962\\
  0.2&2&100&586.125&17.973&0.2008&0.010&0.948\\
  0.1&2&100&658.269&51.069&0.1007&0.006&0.911\\
  0.05&2&100&699.602& 47.911&0.0508&0.004&0.876\\
  0.01&2&100& 1395.3& 367.8& 0.0102&0.001& 0.874\\
  \hline
\end{tabular}
\end{center}
\end{table}

\section*{Acknowledgement}
The authors thank the editor for the thoughtful and constructive comments and suggestions which lead to the improvement of the presentation.

\section*{Appendix}
\appendix
\section{The expected value of $\widehat{p}_{N_2}$}
\label{A:a}
Notice first that
\begin{equation*}
\displaystyle
E\left(\widehat{p}_{N_2}\right)=E\left(E\left(\widehat{p}_{N_2}|
\overline{X}_1\right)\right),
\end{equation*}
where
\begin{equation*}
\displaystyle
E\left(\widehat{p}_{N_2}|\overline{X}_1\right)=\frac{m}{N_2} \widehat{p}_1+\frac{M}{N_2} E\left(\widehat{p}_2|\overline{X}_1\right).
\end{equation*}
Moreover,
\begin{equation*}
\displaystyle
E\left(\widehat{p}_2|\overline{X}_1\right)\approx p+\frac{\left(k_2-1\right)\left(1-\left(1-p\right)^{k_2}\right)}
{2 M \left(k_2\right)^2 \left(1-p\right)^{k_2-1}}.
\end{equation*}
It follows that
\begin{equation*}
\displaystyle
E\left(\widehat{p}_{N_2}\right)=m E\left(\frac{\widehat{p}_1}{N_2}\right)+p E\left(\frac{M}{N_2}\right)+
E\left(\frac{\left(k_2-1\right)\left(1-\left(1-p\right)^{k_2}\right)}
{2 N_2 \left(k_2\right)^2 \left(1-p\right)^{k_2-1}}\right).
\end{equation*}
If $p$ is known, then we can determine the optimal $k$ which is the same for both stages. In this case $k_1=k_2$.
On the other hand, when $p$ is unknown, $k_2$ is a function of $\widehat{p}_1$.
In this case in order to siplify the approximation, we will assume that $k_2$ is a constant independent of $\overline{X}_1$. In this case,
 \begin{equation*}
 \displaystyle
 E\left(\widehat{p}_{N_2}\right)\approx p+m E\left(\frac{\widehat{p}_1}{N_2}\right)+\left(\frac{\left(k_2-1\right)
 \left(1-\left(1-p\right)^{k_2}\right)}{2\left(k_2\right)^2 \left(1-p\right)^{k_2-1}}-m p\right)E\left(\frac{1}{N_2}\right).
 \end{equation*}
 It remains to find an approximation to $E\left(\frac{\widehat{p}_1}{N_2}\right).$
 \begin{equation*}
 \displaystyle
 E\left(\frac{m\widehat{p}_1}{N_2}\right)=\frac{m}{\zeta_{k_1}}
 E\left(\frac{\left(1-\overline{X}_1\right) \left(1-\left(1-\overline{X}_1\right)^{1/k_1}\right)^3}{\overline{X}_1
 \left(1-\overline{X}_1\right)^{2/k_1}}\right).
 \end{equation*}
By the delta method, let
\begin{equation*}
\displaystyle
G\left(x\right)=\frac{\left(1-x\right) \left(1-\left(1-x\right)^{1/k}\right)^3}{x \left(1-x\right)^{2/k}}.
\end{equation*}
Then,
\begin{equation*}
\displaystyle
G^{'}\left(x\right)=\frac{\left(1-\left(1-x\right)^{1/k}\right)^2 \left(\left(1-x\right)^{1/k}\left(x+k\right)-x+2k\right)}{kx^2\left(1-x
\right)^{2/k}}.
\end{equation*}
and
\begin{equation*}
\displaystyle
G^{''}\left(x\right)=\frac{-\left(1-\left(1-x\right)^{1/k}\right)\left(
A\left(x\right)+B\left(x\right) \left(1-x\right)^{1/k}+C\left(x\right)
\left(1-x\right)^{2/k}\right)}{x^3 k^2 \left(1-x\right)^{2/k+1}},
\end{equation*}
where
\begin{equation*}
\displaystyle
A\left(x\right)=-4 x\left(k-x\right)+2 k \left(k\left(1-x\right)+x^2\right)
\end{equation*}
\begin{equation*}
\displaystyle
B\left(x\right)=x k \left(2-x\right)-4 k^2 \left(1-x\right)+x^2
\end{equation*}
\begin{equation*}
\displaystyle
C\left(x\right)=-x^2 \left(k-1\right)+2 k^2 \left(1-x\right)+2 x k
\end{equation*}
Finally, let $g\left(p\right)=G\left(1-\left(1-x\right)^{1/k}\right)$ then
\begin{equation*}
\displaystyle
E\left(\frac{m\widehat{p}_1}{N_2}\right)=\frac{m}{\zeta_1} \left(g\left(p\right)+\frac{1}{2 m} g^{''}\left(p\right) \left(1-p\right)^{k_1} \left(1-\left(1-p\right)^{k_1}\right)\right).
\end{equation*}

\section{The variance of $ \widehat{p}_{N_2}$}
\label{A:b}
The conditional variance of $\widehat{p}_{N_2}$ given $\overline{X}_1$ is
\begin{equation*}
V\left(\widehat{p}_{N_2}|\overline{X}_1\right)=
\left(\frac{M}{N_2}\right)^2 V\left(\widehat{p}_2|\overline{X}_1\right).
\end{equation*}
The conditional expectation of $\widehat{p}_{N_2}$ give $\overline{X}_1$ is
\begin{equation*}
\displaystyle
E\left(\widehat{p}_{N_2}|\overline{X}_1\right)=\frac{m}{N_2} \widehat{p}_1+
\frac{M}{N_2} E\left(\widehat{p}_2|\overline{X}_1\right).
\end{equation*}
We apply the formula $$V\left(\widehat{p}_{N_2}\right)=E\left(V\left(\widehat{p}_{N_2}|
\overline{X}_1\right)\right)+V\left(E\left(\widehat{p}_{N_2}|
\overline{X}_1\right)\right)$$
As derived before, asymptotically
\begin{equation*}
\displaystyle
E\left(\left(\frac{M}{N_2}\right)^2 V\left(\widehat{p}_2|\overline{X}_1\right)\right)\approx\left(\frac{1-p}{k_2}
\right)^2\left(\left(1-p\right)^{-k_2}-1\right) E\left(\frac{M}{N_2^2}\right).
\end{equation*}
The variance of the conditional expectation is more complicated,i.e.
\begin{equation*}
\displaystyle
V\left(\frac{m}{N_2} \widehat{p}_1+\frac{M}{N_2} E\left(\widehat{p}_2|\overline{X}_1\right)\right)=
V\left(\frac{m}{N_2} \widehat{p}_1+\frac{M}{N_2}p+\left(\frac{\left(k_2-1\right) \left(1-\left(1-p\right)^{k_2}\right)}{2\left(k_2\right)^2 \left(1-p\right)^{k_2-1}}\right)^2 \frac{1}{N_2}\right).
\end{equation*}
The right hand term of the last equation is equal to
\begin{align*}
&
V\left(\frac{m \widehat{p}_1}{N_2}\right)+
\left(\left(m p\right)^2-2 p m \frac{\left(k_2-1\right) \left(1-\left(1-p\right)^{k_2}\right)}{2\left(k_2\right)^2 \left(1-p\right)^{k_2-1}}+
\left(\frac{\left(k_2-1\right) \left(1-\left(1-p\right)^{k_2}\right)}{2\left(k_2\right)^2 \left(1-p\right)^{k_2-1}}\right)^2\right) V\left(\frac{1}{N_2}\right)\\
&
-2\left(m^2 p-\frac{\left(k_2-1\right) \left(1-\left(1-p\right)^{k_2}\right)}{2\left(k_2\right)^2 \left(1-p\right)^{k_2-1}}\right)\left(E\left(\frac{\widehat{p}_1}{N_2^2}
\right)-E\left(\frac{\widehat{p}_1}{N_2}\right)E\left(\frac{1}{N_2}\right)
\right).
\end{align*}
According to the previous definition of the function G, and using the delta method we get
\begin{align*}
&
V\left(\frac{m \widehat{p}_1}{N_2}\right)=\frac{m}{\zeta_{k-1}^2}\left(\frac{p}{1-p}
\right)^4\\
&
\frac{\left(1-p\right)^{k_1}\left(\left(1-p\right)\left(k_1+1-\left(1-p
\right)^{k_1}+2\left(1-\left(1-p\right)^{k_1}\right)-k_1\right)\right)}
{k_1^2 \left(1-\left(1-p\right)^{k_1}\right)^3}.
\end{align*}
In order to compute $E\left(\frac{\widehat{p}_1}{N_2^2}\right)$ we introduce the function
\begin{equation*}
\displaystyle
G_2\left(x\right)=\frac{\left(1-x\right)^2 \left(1-\left(1-x\right)^{1/k}\right)^5}{x^2\left(1-x\right)^{4/k}}.
\end{equation*}
The second order derivative of this function is
\begin{equation*}
\displaystyle
G_2^{''}\left(x\right)=-\frac{p^3}{k^2 \left(1-p\right)^4 \left(\theta_k\right)^4}
\left(A\left(\theta_k\right)+B\left(\theta_k\right) \left(1-p\right)+
C\left(\theta_k\right) \left(1-p\right)^2\right).
\end{equation*}
where,
\begin{equation*}
\displaystyle
A\left(\theta_k\right)=4 \theta_k \left(k-\theta_k\right) \left(4+k\right)-6k^2
\end{equation*}
and
\begin{equation*}
\displaystyle
B\left(\theta_k\right)=-4\theta_k\left( k\right) \left(3+2k\right)+3 \theta_k^2\left(k-1\right)+12 k^2
\end{equation*}
and
\begin{equation*}
\displaystyle
C\left(\theta_k\right)=\theta_k\left(k-1\right)\left(4k+\theta_k\right)-6 k^2.
\end{equation*}
Finally,
\begin{equation*}
\displaystyle
E\left(\frac{\widehat{p}_1}{N_2^2}\right)\approx\frac{1}{\zeta_{k_1}}
G_2\left(\theta_{k_1}\right)+\frac{1}{2 m \zeta_k^2}G_2^{''}\left(x\right)\theta_{k_1} \left(1-\theta_{k_1}\right).
\end{equation*}

{}

\end{document}